\documentclass{amsart}

\usepackage{graphicx}
\usepackage{longtable}

\newtheorem{defn}{Definition}

\newtheorem{rem}{Remark}

\def\Li{\mathop{\mathrm{Li}}\nolimits}
\def\Res{\mathop{\mathrm{Res}}\nolimits}

\begin{document}

\title[Integrals of Prime Zeta Function]{Twenty Digits of Some Integrals of the Prime Zeta Function}

\author{Richard J. Mathar}
\urladdr{http://www.mpia-hd.mpg.de/~mathar}
\address{Max-Planck Institute of Astronomy, K\"onigstuhl 17, 69117 Heidelberg, Germany}


\date{\today}

\begin{abstract}
The double sum $\sum_{s\ge 1}\sum_p 1/(p^s \log p^s) = 2.00666645\ldots$
over  the inverse of the product of
prime powers $p^s$ and their logarithms, is computed to 24 decimal digits.
The sum covers all primes $p$ and all integer exponents $s\ge 1$.
The calculational
strategy is adopted from Cohen's work which basically looks at the
fraction as the underivative of the Prime Zeta Function, and then evaluates
the integral by numerical methods.
\end{abstract}

\maketitle
\section{Overview}

\begin{defn}
The constant in the focus of this work is the sum
\begin{equation}
C\equiv \sum_{s=1}^\infty \sum_p \frac{1}{p^s \log(p^s)}
\end{equation}
over all primitive prime powers $p^s$ of the prime numbers $p$
\cite[2.27.3]{Finch}\cite[\S 2.2]{FinchEAMathConst}.
\end{defn}
\begin{rem}
The prime powers are represented by Sloane's sequence A000961---we
drop
the leading term $p^s=1$ to avoid division by zero---.
The constant $C$ and the contribution of $s=1$ are entry A137250
and A137245 \cite{sloane}.
\end{rem}

The only aim of the work is to improve on the previous
estimates
$2.00< C<2.01$
\cite{ErdosCMA26}
and $C>2.006$
\cite{ClarkCMA35}.

The simple computational strategy is to accumulate the 
partial sums
over $s$ in
\begin{equation}
C\equiv \sum_{s=1}^\infty\frac{1}{s}\sum_p \frac{1}{p^s \log p}
,
\label{eq.Csum}
\end{equation}
which converge at reasonable speed
to an accuracy of $10^{-7}$ in $C$ after 20 terms
or to an accuracy of $10^{-19}$ after 59 terms, for example.

\section{Cohen's Integral}
The implementation is a shameless replica of Henri Cohen's reduction of the
double sum over exponents and primes to a series of integrals \cite{Cohen}.
\subsection{Logarithm-to-Integral Conversion}

The logarithm of Euler's formula for Riemann's zeta function
\cite{BorweinJCAM121}\cite[9.523.1]{GR}
\begin{equation}
\zeta(s) = \prod_p \frac{1}{1-p^{-s}}
\label{eq.Eul}
\end{equation}
is
\begin{equation}
\log \zeta(s) = -\sum_p \log (1-p^{-s})
,
\label{eq.logEul}
\end{equation}
which turns with the Taylor expansion of the logarithm into \cite[9.523.2]{GR}\cite{DirichletJM4}
\begin{equation}
\log \zeta(s) = \sum_{k=1}^\infty \sum_p \frac{1}{kp^{ks}}
.
\end{equation}

\begin{defn}
The Prime Zeta Function is \cite{FrobergBIT8,MatharArxiv0803}
\begin{equation}
P(s)=\sum_p\frac{1}{p^s};\quad \Re s>1
.
\label{eq.Pdef}
\end{equation}
\end{defn}

So the penultimate equation can be rephrased as
\begin{equation}
\log \zeta(s) = \sum_{k=1}^\infty \frac{1}{k}P(ks)
.
\end{equation}
The M\"obius inversion of this formula
reads
\cite[(9b)]{BenderAMM82}\cite[\S 17.1.3]{ErdelyiIII}\cite{SebahGourdon}
\begin{equation}
P(s)=\sum_{k=1}^\infty \frac{\mu(k)}{k} \log \zeta(ks).
\label{eq.mob}
\end{equation}
First integrals of terms in (\ref{eq.Pdef}) are \cite[565.1.]{Dwight}
\begin{equation}
\int_s^\infty \frac{1}{p^t}dt = \frac{1}{p^s \log p}
,\quad s>1,
\end{equation}
so integration of (\ref{eq.mob}) using (\ref{eq.Pdef}) gives
\begin{equation}
\int_s^\infty \sum_p\frac{dt}{p^t}
=
\sum_p\frac{1}{p^s\log p}
=
\sum_{k=1}^\infty \frac{\mu(k)}{k} \int_s^\infty \log \zeta(kt)dt
=
\sum_{k=1}^\infty \frac{\mu(k)}{k^2} \int_{ks}^\infty \log \zeta(t)dt
.
\label{eq.logs}
\end{equation}
\begin{defn}
Cohen's integral is
\begin{equation}
I_n(m)\equiv \frac{1}{n!} \int_m^\infty (t-m)^n \log \zeta(t)dt,\quad n=0,1,2,\ldots
\end{equation}
over the logarithm of Riemann's zeta function with variable integer lower limit \cite{Cohen}.
\end{defn}
Insertion into (\ref{eq.logs}) outlines the strategy to evaluate (\ref{eq.Csum}),
\begin{equation}
\sum_p\frac{1}{p^s\log p}
=\sum_{k=1}^\infty \frac{\mu(k)}{k^2}I_0(ks)
.
\label{eq.logsum}
\end{equation}

\subsection{Transition to the Dirichlet Eta Function}
For $I_0(1)$, the pole of $\zeta(t)$ at $t=1$ can be isolated by handling the singularity
via the smooth function \cite[23.2.19]{AS}\cite[p 267]{WhittakerW}
\begin{equation}
\eta(s)\equiv (1-2^{1-s})\zeta(s)
.
\label{eq.etaDef}
\end{equation}
For lower limits of the integral larger than one, this recipe is not needed;
we stick to it to present a shorter, simpler program.
\begin{rem}
The finite value at $s=1$ is \cite{Finch}\cite[0.232.1,1.511]{GR}
\begin{equation}
\eta(1)=\log 2.
\end{equation}
This follows also from the residuum of $\Res \zeta(s)_{\mid s=1}=1$
in conjunction with the Taylor expansion
$1-2^{1-s}= -\sum_{l=1}^\infty \frac{(-\log 2)^l}{l!}(s-1)^l=(s-1)\log 2-\frac{\log^2 2}{2}(s-1)^2+O((s-1)^3)$.
If we introduce
Stieltjes constants $\gamma_j$ \cite{CoffeyJMA317,CoffeyPRSA462,Coffeyarxiv07},
\begin{equation}
\zeta(s)=\frac{1}{s-1}+\gamma-\gamma_1(s-1)+\frac{\gamma_2}{2}(s-1)^2+\cdots,
\end{equation}
the Taylor series of $\eta$ near $s=1$ becomes
\begin{equation}
\eta(s)=\log 2+\log 2\left(\gamma-\frac{\log 2}{2}\right)(s-1)
-\log 2\left(\gamma_1+\gamma\frac{\log 2}{2}-\frac{\log^2 2}{6}\right)(s-1)^2
+\ldots
.
\end{equation}
\end{rem}

The integrated logarithm of (\ref{eq.etaDef}) is
\begin{equation}
\int_m^\infty ds \log \eta(s)
=
\int_m^\infty ds \log (1-2^{1-s})
+
\int_m^\infty ds \log \zeta(s)
.
\label{eq.logsplit}
\end{equation}
One term is a dilogarithm \cite{GinsbergCACM18,KirillovPTP118,KolbigBIT10,LoxtonAA43},
\begin{equation}
\int_m^\infty \log\left(1-2^{1-s}\right)ds
=
\int_{1-2^{1-m}}^1 \frac{\log x}{\log 2(1-x)}dx
\\
=
- \frac{1}{\log 2}\Li_2\left(\frac{1}{2^{m-1}}\right) .
\end{equation}

\begin{rem}
Special values at $m=1$ and $m=2$ are $\Li_2(1)=\pi^2/6$
and $\Li_2(1/2)=\pi^2/12-\frac{1}{2}\log^2(2)$ \cite{Lewin};
see the constants A013661 and A076788 in the On-Line Encyclopedia of
Integer Sequences \cite{sloane}.
An accurate representation of
$\pi^2/(12\log 2)$ is entry
A100199.
$\eta(1)=\log 2$ is A002162, $\gamma$ is A001620,
$\eta'(1)=\log 2(\gamma-\log 2/2)$ is A091812,
and $-\gamma_1$ is A082633.
\end{rem}
This moulds (\ref{eq.logsplit}) into
\begin{equation}
I_0(m)
=
\int_m^\infty ds \log \eta(s)
+ \frac{1}{\log 2}\Li_2\left(\frac{1}{2^{m-1}}\right)
.
\label{eq.Isplit}
\end{equation}

\subsection{Numerical Implementation}
The variable substitution
\begin{equation}
u=1-\frac{m}{s};\quad s=\frac{m}{1-u};\quad ds=\frac{m}{(1-u)^2}du
\label{eq.usubs}
\end{equation}
maps the interval $m\le s<\infty$ onto the interval $0\le u\le 1$.
\begin{rem}
Alternative substitutions like $u=1/m-m/s^2$ or $u=1-2^{-s}$ also map
the half-infinite $s$-interval to a finite $u$-interval. They have the advantage
of more evenly balanced integral kernels, but the disadvantage of
infinite slopes at one end of the $u$-interval.
\end{rem}
Wynn's $\epsilon$-algorithm
\cite{Wynn1956,LevinAMC9}
is applied to numerical values gathered by
the trapezoidal rule to evaluate the $\eta$-integral in (\ref{eq.Isplit}),
\begin{equation}
\int_m^\infty ds \log \eta(s)
=
m\int_0^1 \frac{ \log\eta\left(\frac{m}{1-u}\right)}{(1-u)^2}du
.
\label{eq.Ieta}
\end{equation}

\begin{figure}
\includegraphics[scale=0.7]{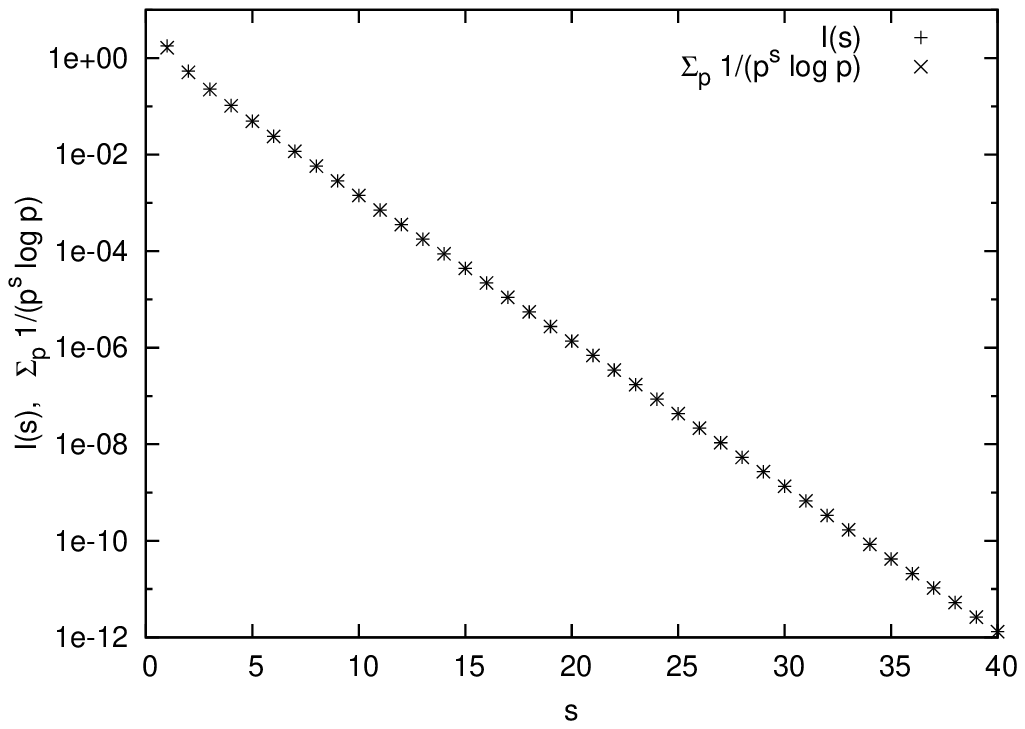}
\caption{
The functions $I_0(s)$ and $\sum_p \frac{1}{p^s\log p}$
fall of exponentially as a function of their argument $s$,
and are not distinguishable in this common plot.
}
\label{fig.Itab}
\end{figure}

\begin{rem}
Because $\log \zeta(s)\approx 2^{-s}$ at $s>10$, $I_0(m)\approx 1/(2^m \log2)$ as $m\to \infty$.
This functional dependence leads to an almost straight line on the semi-logarithmic plot
in Fig.\ \ref{fig.Itab}. So the first neglected term in (\ref{eq.Csum}) is a good
estimator to the error in the partial sums. 
\end{rem}

\begin{figure}
\includegraphics[scale=0.7]{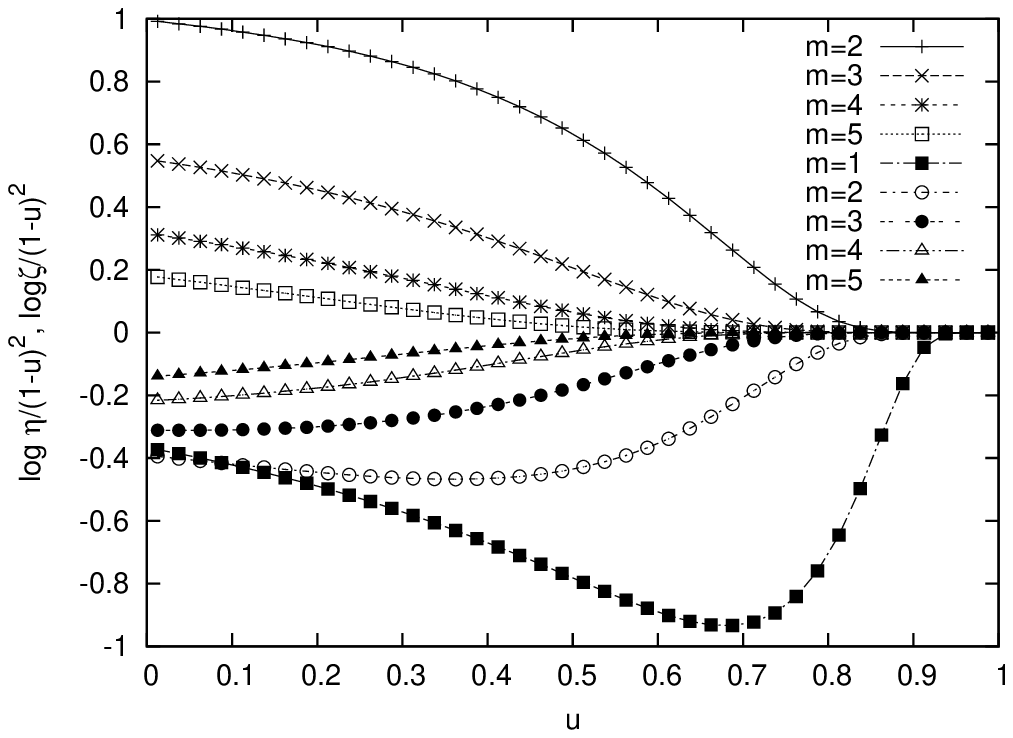}
\caption{
The functions $( m\log \eta)/(1-u)^2 <0$ and
$(m\log \zeta)/(1-u)^2>0$
arising in the integral kernels through the substitution (\ref{eq.usubs}).
}
\end{figure}

\subsection{Intermediate Results}
The following table shows $s$ and $I_0(s)$ in the left double column,
$s$ and $\sum_p 1/(p^s\log p)$ in the right double column, last digits rounded:
\small \begin{verbatim}
  1 1.797569958628739407930251e+00    1 1.636616323351260868569658e+00
  2 5.365269459211771009617190e-01    2 5.077821878591993187743751e-01
  3 2.275044530363196437654595e-01    3 2.212033403969814969599433e-01
  4 1.041502253168599126451383e-01    4 1.026654782331571962664758e-01
  5 4.942465209576037094607771e-02    5 4.906357474830439568004180e-02
  6 2.392429939526549813915300e-02    6 2.383519825459824737569112e-02
  7 1.171800374394507398697301e-02    7 1.169586557789515115400634e-02
  8 5.781462903507255601662420e-03    8 5.775944595316604352784204e-03
  9 2.865717415930253367688719e-03    9 2.864339771659922115926198e-03
 10 1.424706500149304604154559e-03   10 1.424362320209387639441021e-03
 11 7.096782687818267364308355e-04   11 7.095922515383712411613441e-04
 12 3.539573278954885950994965e-04   12 3.539358269254206683875757e-04
 13 1.766870488681893665136623e-04   13 1.766816740292220402778692e-04
 14 8.824687187128450251937842e-05   14 8.824552821042831521626792e-05
 15 4.409135113515804839684079e-05   15 4.409101522588177584140913e-05
 16 2.203501275340515891846297e-05   16 2.203492877680910955374197e-05
 17 1.101395992552518507579192e-05   17 1.101393893146441713313451e-05
 18 5.505799725855358158700238e-06   18 5.505794477350959773416061e-06
 19 2.752506920732716039062161e-06   19 2.752505608607939317868454e-06
 20 1.376122595515799546777009e-06   20 1.376122267484767801052807e-06
 21 6.880177047831559600925881e-07   21 6.880176227754180030145763e-07
 22 3.439943284950306802000172e-07   22 3.439943079930986522426868e-07
 23 1.719923247093458602614288e-07   23 1.719923195838631569242761e-07
 24 8.599454961284003932811501e-08   24 8.599454833146940100459182e-08
 25 4.299673733467609442174442e-08   25 4.299673701433343948018079e-08
 26 2.149818953720244919068455e-08   26 2.149818945711678602970514e-08
 27 1.074903506531325106214723e-08   27 1.074903504529183534309235e-08
 28 5.374497633245768258586266e-09   28 5.374497628240414337653354e-09
 29 2.687242183906380835200363e-09   29 2.687242182655042356063437e-09
 30 1.343618881152634756065114e-09   30 1.343618880839800136417081e-09
 31 6.718087036690018174171263e-10   31 6.718087035907931625220498e-10
 32 3.359041062052731307704559e-10   32 3.359041061857209670487929e-10
 33 1.679519712278619904984955e-10   33 1.679519712229739495683419e-10
 34 8.397595832274787450826827e-11   34 8.397595832152586427576250e-11
 35 4.198797006441494468678066e-11   35 4.198797006410944212865828e-11
 36 2.099398199991330939924248e-11   36 2.099398199983693375971239e-11
 37 1.049698998919830689452100e-11   37 1.049698998917921298463855e-11
 38 5.248494657681297045976310e-12   38 5.248494657676523568505703e-12
 39 2.624247216535094520600355e-12   39 2.624247216533901151232704e-12
 40 1.312123570832462108601650e-12   40 1.312123570832163766259737e-12
 41 6.560617729378942125931467e-13   41 6.560617729378196270076685e-13
 42 3.280308823095077105921359e-13   42 3.280308823094890641957663e-13
 43 1.640154397682756111408190e-13   43 1.640154397682709495417266e-13
 44 8.200771942197877942789073e-14   44 8.200771942197761402811764e-14
 45 4.100385955693647813828150e-14   45 4.100385955693618678833822e-14
 46 2.050192972711729282745366e-14   46 2.050192972711721998996784e-14
 47 1.025096484644167040382940e-14   47 1.025096484644165219445795e-14
 48 5.125482417515178049571403e-15   48 5.125482417515173497228539e-15
 49 2.562741206855703686734819e-15   49 2.562741206855702548649103e-15
 50 1.281370602793890158864587e-15   50 1.281370602793889874343158e-15
 51 6.406853011856245416428711e-16   51 6.406853011856244705125138e-16
 52 3.203426505223720974861043e-16   52 3.203426505223720797035150e-16
 53 1.601713252377059924465468e-16   53 1.601713252377059880008995e-16
 54 8.008566261102631116158675e-17   54 8.008566261102631005017492e-17
 55 4.004283130290426065285105e-17   55 4.004283130290426037499809e-17
 56 2.002141565058249870693295e-17   56 2.002141565058249863746971e-17
 57 1.001070782500137215275764e-17   57 1.001070782500137213539183e-17
 58 5.005353912404060344256377e-18   58 5.005353912404060339914924e-18
 59 2.502676956169821595115832e-18   59 2.502676956169821594030469e-18
 60 1.251338478074174605310912e-18   60 1.251338478074174605039571e-18
 61 6.256692390335085719290662e-19   61 6.256692390335085718612310e-19
 62 3.128346195155613757280562e-19   62 3.128346195155613757110974e-19
 63 1.564173097573830511199490e-19   63 1.564173097573830511157093e-19
 64 7.820865487855897997896814e-20   64 7.820865487855897997790822e-20
 65 3.910432743923530812923694e-20   65 3.910432743923530812897196e-20
 66 1.955216371960292677789151e-20   66 1.955216371960292677782526e-20
 67 9.776081859796554293375620e-21   67 9.776081859796554293359059e-21
 68 4.888040929896640781499146e-21   68 4.888040929896640781495005e-21
 69 2.444020464947774935687030e-21   69 2.444020464947774935685995e-21
 70 1.222010232473705649489420e-21   70 1.222010232473705649489161e-21
 71 6.110051162367922186267000e-22   71 6.110051162367922186266353e-22
 72 3.055025581183759072740187e-22   72 3.055025581183759072740025e-22
 73 1.527512790591812196239003e-22   73 1.527512790591812196238962e-22
 74 7.637563952958836514091412e-23   74 7.637563952958836514091310e-23
 75 3.818781976479343434677847e-23   75 3.818781976479343434677822e-23
 76 1.909390988239646776549639e-23   76 1.909390988239646776549633e-23
 77 9.546954941198150746783920e-24   77 9.546954941198150746783904e-24
 78 4.773477470599047661403870e-24   78 4.773477470599047661403866e-24
 79 2.386738735299514593372572e-24   79 2.386738735299514593372571e-24
 80 1.193369367649754217576498e-24   80 1.193369367649754217576498e-24
\end{verbatim}\normalsize
At large $s$,
the two values equalize because the first term at $k=1$ dominates
the sum at the right hand side of (\ref{eq.logsum}).
Cohen reported the first two rows \cite{Cohen}.
Erd\H{o}s and Zhang published an
upper bound of $\sum_p 1/(p\log p)$ \cite{ErdosAMS117};
a proof of convergence had been given earlier \cite{ErdosPLMS10}.

The rightmost column $P^L(s,1)$ is an aid to computation of the more general
\begin{equation}
P^L(s,a)\equiv \sum_p \frac{1}{(a-1+p)^s\log p}
=
\sum_{l=0}^\infty \binom{s+l-1}{s-1}(1-a)^l P^L(l+s,1).
\end{equation}
The $l$-series converges for $a=0$ and $a=2$.
Larger $a$ are then reached recursively
with the alternating geometric series
\begin{equation}
P^L(s,a+1)
=\sum_p \frac{1}{(a-1+p)^s\left(1+\frac{1}{a-1+p}\right)^s\log p}
=
\sum_{l=0}^\infty \binom{s+l-1}{s-1}(-1)^l P^L(s+l,a).
\end{equation}
This yields the following short table with double
columns of $s$, $a$ and $P^L(s,a)$ each:
\small \begin{verbatim}
 1  0 2.564343220686309193e+00       4  2 2.200937333087681441e-02
 2  0 1.735535173734295407e+00       5  2 6.924273989916642133e-03
 3  0 1.569430040669505312e+00       6  2 2.216719348246691407e-03
 4  0 1.502480490856796843e+00       7  2 7.177059901294852453e-04
 5  0 1.471819234100400691e+00       8  2 2.341804937450834077e-04
 6  0 1.457080819283874067e+00       9  2 7.683437427098209717e-05
 7  0 1.449846098516752367e+00      10  2 2.531100154045370884e-05
 8  0 1.446260454845751749e+00      11  2 8.362846790350413196e-06
 9  0 1.444475273576531254e+00      12  2 2.769232699593665178e-06
10  0 1.443584547482253809e+00      13  2 9.185072149314254133e-07
11  0 1.443139643195342866e+00      14  2 3.050306451883795961e-07
12  0 1.442917304533643634e+00      15  2 1.013929603798279887e-07
13  0 1.442806163373811902e+00      16  2 3.372678630308857437e-08
14  0 1.442750599803600762e+00      17  2 1.122456404716334707e-08
15  0 1.442722819765431096e+00      18  2 3.737099942191877583e-09
16  0 1.442708930182166938e+00      19  2 1.244595144792613520e-09
17  0 1.442701985499337981e+00      20  2 4.145889249942980892e-10
18  0 1.442698513185098959e+00       1  3 1.083884561790478549e+00
19  0 1.442696777034769093e+00       2  3 1.539997501094314544e-01
20  0 1.442695908961300868e+00       3  3 3.279427700685188571e-02
 1  2 1.280433764453964036e+00       4  3 7.457583410785408360e-03
 2  2 2.518771844123868571e-01       5  3 1.747751882461897344e-03
 3  2 7.208846839108258174e-02
\end{verbatim}\normalsize
In the left column of this table we observe $\lim_{s\to \infty} P^L(s,0)=1/\log(2)\approx 1.44\ldots$.

\section{Higher Powers of The Logarithms}
\subsection{Quadratic}

Repeated integration in the spirit of \eqref{eq.Csum} generates higher
powers of the logarithm in the denominator of the prime sums:

\begin{equation}
\int_s^\infty \frac{1}{p^u \log p}du = \frac{1}{p^s \log^2 p}
= \int_s^\infty du \int_u ^\infty dt \frac{1}{p^t},
\end{equation}
which integrates an infinite triangular region which extends right from $t=u$
and upwards from the diagonal $u\ge s$.
Interchange of the order of the two integrations contracts it
to a single integral:
\begin{multline}
\sum_p \frac{1}{p^s \log^2 p}
= \int_s^\infty du \int_u ^\infty dt P(t)
= \int_{t=s}^\infty dt \int_{u=s}^t ds P(t)\\
= \int_{t=s}^\infty dt (t-s) \sum_{k=1}^\infty \frac{\mu(k)}{k}\log \zeta(kt)
= \int_{x=sk}^\infty dx (x-sk) \sum_{k=1}^\infty \frac{\mu(k)}{k^3}\log \zeta(x) \\
= 
\sum_{k=1}^\infty \frac{\mu(k)}{k^3} I_1(sk)
.
\label{eq.log2}
\end{multline}
As before it is numerically advantageous to employ the $\eta$-function \eqref{eq.etaDef}
in the integral kernel of $I_1$ for small $m$:
\begin{equation}
\int_m^\infty dt (t-m)\log \eta(t)
=
\int_m^\infty dt (t-m)\log (1-2^{1-t})
+
\int_m^\infty dt (t-m)\log \zeta(t).
\end{equation}
The contribution from $1-2^{1-t}$ becomes a mix of well-converging polylogarithmic constants.
We substitute
\begin{equation}
x=1-2^{1-t},\quad t=1-\log(1-x)/\log 2,\quad dt=\frac{dx}{(1-x)\log 2}
\label{eq.xlog2}
\end{equation}
to approach standard integrals over products of logarithms:
\begin{multline}
\int_m^\infty dt (t-m)\log (1-2^{1-t})
\\
=
\frac{1}{\log 2}\int_{1-2^{1-m}}^1 \frac{dx}{1-x} [1-\frac{\log(1-x)}{\log 2}-m]\log x\\
=
\frac{1-m}{\log 2}\int_0^{2^{1-m}}\frac{du}{u}\log(1-u)
-
\frac{1}{\log^2 2}
\int_0^{2^{1-m}} \frac{du}{u} \log(u)\log (1-u)\\
=
\frac{1-m}{\log 2} l_{0,1}(2^{1-m})
-
\frac{1}{\log^2 2} l_{1,1}(2^{1-m}).
\end{multline}
The $l$-integrals are defined by the powers of the two logarithms
in the kernel, and discussed in the appendix for the upper limits of interest:
\begin{defn} (Polylogarithmic Constants)
\begin{equation}
l_{i,j}(x) \equiv \int_0^x \log^i (u) \log^j (1-u)\frac{du}{u}.
\label{eq.ldef}
\end{equation}
\end{defn}

The values of $I_1$ in \eqref{eq.log2} are eventually computed
for $m=1$ and $m=2$ via
\begin{equation}
I_1(m) = \int_m^\infty (s-m)\log \eta(s)ds -\int_m^\infty (s-m)\log(1-2^{1-s})ds
\label{eq.I1defer}
\end{equation}
and numerical integration of the first term on the right hand side.
The important results are those for small $s$ where direct summation over the primes
suffers from slow convergence
\cite[A319231,A319232]{sloane}:
\begin{eqnarray}
\sum_p \frac{1}{p\log^2 p}&\approx& 1.520970439939500\ldots \\
\sum_p \frac{1}{p^2\log^2 p}&\approx& 0.637056184074676 \ldots\\
\sum_p \frac{1}{p^3\log^2 p}&\approx& 0.2949791653099842 \ldots
\end{eqnarray}

\subsection{Cubic}
Integrating \eqref{eq.log2} once more lifts the square to a cube in the denominator:
\begin{multline}
\sum_p \frac{1}{p^s \log^3 p}
= \int_s^\infty du \sum_p \frac{1}{p^u \log^2 p}\\
= \int_s^\infty du \int_{t=u}^\infty dt (t-u)\sum_{k=1}^\infty \frac{\mu(k)}{k}\log \zeta(kt)\\
= 
\int_{t=s}^\infty dt \int_{u=s}^t du (t-u)\sum_{k=1}^\infty \frac{\mu(k)}{k}\log \zeta(kt)\\
= 
\frac12 \int_{s}^\infty dt (t-s)^2\sum_{k=1}^\infty \frac{\mu(k)}{k}\log \zeta(kt)
= \frac12 \int_{ks}^\infty dx (x-ks)^2 \sum_{k=1}^\infty \frac{\mu(k)}{k^4}\log \zeta(x) \\
= \sum_{k=1}^\infty \frac{\mu(k)}{k^4} I_2(sk).
\end{multline}
\begin{equation}
\sum_p \frac{1}{p^s\log^n p} =\sum_{k=1}^\infty \frac{\mu(k)}{k^{n+1}}I_{n-1}(sk),\quad n=1,2,3,\ldots .
\end{equation}
Similar to \eqref{eq.I1defer} we defer the $\zeta$-function to the $\eta$-function for small $m$:
\begin{equation}
I_2(m)= \frac12 \int_{s=m}^\infty (s-m)^2\log \eta(s) ds
-
\frac12 \int_{s=m}^\infty (s-m)^2\log (1-2^{2-m}) ds.
\end{equation}
The second term is sliced with \eqref{eq.xlog2} as
\begin{multline}
\frac12 \int_{s=m}^\infty (s-m)^2\log (1-2^{1-m}) ds
=
\frac{1}{2\log 2} \int_{x=1-2^{1-m}}^1 \left[1-\frac{\log(1-x)}{\log 2}-m\right]^2\log x \frac{dx}{(1-x)}
\\
=
\frac{1}{2\log 2} \int_{u=0}^{2^{1-m}} [1-\frac{\log(u)}{\log 2}-m]^2\log (1-u) \frac{du}{u}
\\
=
\frac{(1-m)^2}{2\log 2} \int_{u=0}^{2^{1-m}} \log (1-u) \frac{du}{u}
-
\frac{1-m}{\log^2 2} \int_{u=0}^{2^{1-m}} \log u \log (1-u) \frac{du}{u}
+
\frac{1}{2\log^3 2} \int_{u=0}^{2^{1-m}} \log^2 u\log (1-u) \frac{du}{u}
\\
=
\frac{(1-m)^2}{2\log 2} l_{0,1}(2^{1-m})
-
\frac{1-m}{\log^2 2} l_{1,1}(2^{1-m})
+
\frac{1}{2\log^3 2} l_{2,1}(2^{1-m})
,
\end{multline}
such that for $m=1$ and $m=2$ the $l_{i,j}$ can be taken from the appendix.
\begin{eqnarray}
\sum_p \frac{1}{p\log^3 p}&\approx& 1.8461474193664495 \ldots\\
\sum_p \frac{1}{p^2\log^3 p}&\approx& 0.848270491703549\ldots\\
\sum_p \frac{1}{p^3\log^3 p}&\approx& 0.405696596756787 \ldots
\end{eqnarray}

\section{Summary}
We obtained the constant
\begin{equation}
 C= 2.00666645283106875643229\ldots
\end{equation}
by a rather basic numerical approach to Cohen's integral over the logarithm of the Prime Zeta Function.
A table of $\sum_p 1/(p^s \log p)$ was generated for $1\le s\le 80$
and a table of $\sum_p 1/[(p\pm 1)^s \log p]$ for $1\le s\le 20$.

\begin{appendix}
\section{Polylogarithm Constants}
We gather the values of $l_{i,j}(x)$ that emerge
in the main text while evaluating $I_j$ for $j=0,1$ and 2 and $m=1$ or 2.
If $j=0$ the integrals are elementary \cite[616.1]{Dwight}:
\begin{equation}
\int \log ^i(u) \frac{du}{u} = \frac{\log^{i+1}(u)}{i+1}.
\end{equation}
If $j=x=1$, the $l_{i,1}(1)$ are essentially $\zeta(i+2)$-values:
\cite{RutledgeAMM45,KolbigMathComp39}
\begin{equation}
l_{1,1}(1)
=
\int_0^1 \frac{du}{u}\log(u)\log(1-u)
=\zeta(3);
\end{equation}
\begin{equation}
l_{2,1}(1) = \int_0^1 \frac{du}{u} \log^2(u) \log(1-u) = -2\zeta(4) = -\frac{\pi^4}{45}.
\end{equation}
If $x=1/2$ we have
\cite[4.291.3]{GR}
\begin{equation}
l_{0,1}(1/2)
=
\int_0^{\frac12} \frac{du}{u}\log(1-u)
=\frac12 \log^2(2)-\frac{\pi^2}{12}.
\end{equation}
A particular case of \cite[A.3.5(8)]{Lewin}
\begin{equation}
l_{1,1}(x) 
= -\log(x)\Li_2(x)+\Li_3(x).
\end{equation}
in conjunction with \cite[A.2.6.(3)]{Lewin}\cite[A099217]{sloane}
\begin{equation}
\Li_3(1/2) = \frac78 \zeta(3)-\frac{\pi^2}{12}\log 2+\frac16 \log^3(2)
\approx 0.537213193608\ldots
\end{equation}
is \cite[A2.1.(4), A2.6(3)]{Lewin}
\begin{equation}
l_{1,1}(1/2)=
-\frac{1}{3} \log^3(2)
+\frac78 \zeta(3)
\approx 0.9407915729350018398184\ldots
\end{equation}

The indefinite integral \cite[A3.5.(10)]{Lewin}
\begin{equation}
l_{2,1}(x) = -2\Li_4(x)+2\log(x)\Li_3(x)-\log^2(x)\Li_2(x)
\end{equation}
generates
\begin{multline}
l_{2,1}(1/2) = -2\Li_4(1/2)-2\log(2)\Li_3(1/2)-\log^2(2)\Li_2(1/2)\\
\approx -2.059432960131268548157225796\ldots
\end{multline}
where \cite[A099218]{sloane}
\begin{equation}
2\Li_4(1/2) = {}_5F_4(1,1,1,1,1;2,2,2,2;1/2) \approx 1.034958123347798772\ldots
\end{equation}

\end{appendix}

\bibliographystyle{amsplain}
\bibliography{all}

\end{document}